\documentclass[a4paper,11pt,twoside,reqno]{amsart}
\usepackage[dvips]{graphicx}
\usepackage[utf8]{inputenc}
\usepackage[plainpages=false,pdfpagelabels=true]{hyperref}
\usepackage{amssymb,amsthm}
\usepackage[margin=1in]{geometry}

   \newtheorem{theorem}{Theorem}[section]
   \newtheorem{proposition}[theorem]{Proposition}

 \theoremstyle{remark}
   
   \newtheorem{remark}[theorem]{Remark}
   
\numberwithin{equation}{section}
\numberwithin{figure}{section}

\allowdisplaybreaks[1]

\numberwithin{equation}{section}

\title{Magnetic geodesics on surfaces with singularities}
\author{Volker Branding and Wayne Rossman}

\date{\today}
\thanks{This research was supported the joint Austrian-Japanese grant \emph{I1671-N26: Transformations and Singularities}.}
\address{TU Wien\\
Institut für diskrete Mathematik und Geometrie\\
Wiedner Hauptstraße 8–10, A-1040 Wien}
\email[]{volker@geometrie.tuwien.ac.at}
\address{Department of Mathematics, Faculty of Science, University of Kobe\\
Rokko, Kobe, Japan 657-8501}
\email[]{wayne@math.kobe-u.ac.jp}

\begin{document}

\begin{abstract}
We prove that, generically, magnetic geodesics on surfaces will turn away from 
points with lightlike tangent planes, 
and we motivate our result with numerical solutions for closed magnetic geodesics.  
\end{abstract} 

\maketitle

\section{Introduction}

A magnetic geodesic describes the trajectory of a charged particle in a 
Riemannian manifold $M$ under the influence of an external magnetic field.
Numerical experimentation suggests that almost all magnetic geodesics tend to 
avoid any lightlike singularities (points where the tangent spaces are lightlike) 
that $M$ may have, regardless of choice of bounded 
smooth external magnetic field.  Our primary result is a mathematically 
rigorous confirmation of this behavior.  

Initially, we take $M$ to be a complete, orientable Riemannian manifold without boundary 
of dimension $n$ and Riemannian metric $\langle \cdot , \cdot \rangle$.  
For a given two form $\Omega$ defined on $M$ we associate a smooth section \(Z \in \text{Hom}(TM,TM)\) defined via
\[
\langle\eta,Z(\xi)\rangle=\Omega(\eta,\xi)
\]
for all $\eta,\xi\in TM$.
We will investigate the existence of closed curves 
$\gamma=\gamma(t)$ satisfying the following equation
\begin{equation}
\label{magnetic-geodesic}
\nabla_{\gamma^\prime} \gamma^\prime = Z(\gamma^\prime).
\end{equation}
Note that, in contrast to geodesics, which correspond to \(Z=0\), the equation for magnetic
geodesics is not invariant under rescaling of \(t\).

In the case that \(M\) is a surface, that is $n=2$, we know that every 
two-form $\Omega$ is a multiple of the volume form $\Omega_0$ associated with 
$\langle \cdot , \cdot \rangle$. Thus, every 
two-form can be written as $\Omega=\kappa\Omega_0$ for some function 
$\kappa\colon M\to\mathbb{R}$. We can exploit this fact
to rewrite the right hand side of \eqref{magnetic-geodesic} as
\begin{equation}
\label{equation-surface}
Z(\gamma')=\kappa J^{90}_\gamma(\gamma'), 
\end{equation}
where $J^{90}_{\gamma}$ represents rotation in the tangent space $T_\gamma M$ 
by angle $\pi/2$, see \cite{Miranda}. 
Due to this fact one often refers to \eqref{equation-surface} as the 
\emph{prescribed geodesic curvature equation}, and $\kappa$ is 
proportional to the geodesic curvature function.  

We will always assume that $\kappa$ is a smooth and bounded function.

\begin{remark}
If the two-form \(\Omega\) is exact, then \eqref{magnetic-geodesic}
also arises from a variational principle, see \cite{cmp}, \cite{Tai2}.
\end{remark}

Note that a solution of \eqref{magnetic-geodesic} has constant speed, which follows from
\begin{equation}\label{eqn1pt3}
\frac{\partial}{\partial t}\frac{1}{2}|\gamma'|^2=\langle\nabla_{\gamma^\prime}\gamma^\prime,\gamma^\prime\rangle
=\langle Z(\gamma'),\gamma^\prime\rangle=\Omega(\gamma^\prime,\gamma^\prime)=0
\end{equation}
due to the skew-symmetry of the two-form \(\Omega\).

For magnetic geodesics on surfaces, several existence results are available, 
employing techniques from symplectic geometry \cite{Gin1}, \cite{Gin2}
and from the calculus of variations \cite{Tai2}.
In the papers of Schneider \cite{SchS2}, \cite{SchH2}, and the paper by Schneider and Rosenberg \cite{SchRos}, 
existence results for closed magnetic geodesics on Riemann surfaces are 
given by studying the zeros of a certain vector field.

Here rather, we give an approach more aimed at usefulness for numerics, and then 
proceed to produce examples of closed magnetic geodesics numerically.  We then   
study the behavior of magnetic geodesics near singular points of a surface by proving 
our main result Theorem \ref{thm:main}, and our proof employs the fact that 
magnetic geodesics have constant speed parameterizations.  

This article is organized as follows:
In section 2 we derive several numerical examples of magnetic geodesics.
Moreover, we provide several analytic statements that support our numerical calculations.
In section 3 we focus on magnetic geodesics on almost-everywhere-spacelike 
surfaces with lightlike 
singularities and show that they will tend to turn away from the singularities 
unless they enter the singular sets at specific angles, which is the content of 
Theorem \ref{thm:main}.  

\section{Closed magnetic geodesics on surfaces in 
Euclidean and Minkowski \texorpdfstring{$3$}{}-spaces}
Before we turn to the numerical integration of \eqref{magnetic-geodesic} let us make the
following observations.

By the Theorem of Picard-Lindeloef we always get a local solution to
\eqref{magnetic-geodesic}. However, similar to the classical Hopf-Rinow
theorem in Riemannian geometry we can show

\begin{theorem}\label{thm:1}
Let \( (M,\langle\cdot,\cdot\rangle) \) be a complete Riemannian surface and
$\kappa : M \to \mathbb{R}$ be a prescribed function.  
Let \[ \gamma(t) : (a,b) \to M \] be a curve in $M$ with geodesic 
curvature $\kappa (\gamma(t) )$ at $\gamma(t)$, in other words, 
$\gamma$ is a nontrivial solution to 
\begin{equation}\label{eqn:1} 
\nabla_{\gamma^\prime} \gamma^\prime = \kappa J^{90}_\gamma(\gamma') \; . 
\end{equation}
Then the domain $(a,b)$ can be extended to all of $\mathbb{R}$.  
\end{theorem}

\begin{proof}
To show that the maximal interval of existence of \eqref{eqn:1} is indeed all 
of \(\mathbb{R}\) we assume that there is a maximal interval of existence 
and then show that we can extend the solution beyond that interval.
Thus, assume that \(\gamma\colon (a,b)\to M\) is a magnetic geodesic 
with maximal domain of definition. Since \(|\gamma'|^2\) is constant we know that the curve \(\gamma\) has constant length \(L[\gamma]\).
Then we have for a sequence \(\gamma(t_i)_{i\in\mathbb{N}}\)
\[
d(\gamma(t_i),\gamma(t_j))\leq L[\gamma_{[t_i,t_j]}]\leq C|t_i-t_j|,
\]
where \(d\) denotes the Riemannian distance function. Hence, \(\gamma(t_i)_{i\in\mathbb{N}}\) is a Cauchy sequence with respect to \(d\).
It is easy to see that the limit is independent of the chosen sequence.

As a next step, we show that we may extend \(\gamma'\) to \((a,b]\). 
To this end we use the local expression for \eqref{eqn:1}, that is
\[
(\gamma^{\prime\prime})^k=-\sum_{i,j=1}^2\Gamma^k_{ij}(\gamma^\prime)^i(\gamma^\prime)^j-\kappa (J^{90}_\gamma(\gamma'))^k,\qquad k=1,2.
\]
Now, consider the expression
\begin{equation*}
|\gamma^\prime(t_i)-\gamma^\prime(t_j)|_{L^\infty}=
\big|\int_{t_i}^{t_j}\gamma^{\prime\prime}(\tau)d\tau\big|_{L^\infty} 
\leq C|t_i-t_j|_{L^\infty}\,.
\end{equation*}
Using that \(|\gamma^\prime|\) is constant it follows 
that \(\gamma^\prime(t_i)\) forms a Cauchy sequence and 
converges to some \(\gamma^\prime_\infty\).
Again, the limit is independent of the chosen sequence.

By differentiating the equation for magnetic geodesics and using the same method as for estimating $|\gamma^\prime(t_i)-\gamma^\prime(t_j)|_{L^\infty}$
we can show that also \(\gamma^{\prime\prime}(t_i)\) forms a Cauchy sequence.

Now, assume that \(\tilde{\gamma}\colon(\beta-a,\beta+a)\to M\) is a magnetic 
geodesic with \(\tilde{\gamma}(\beta)=\hat{\gamma}(\beta)\)
and \(\tilde{\gamma}^\prime(\beta)=\hat{\gamma}^\prime(\beta)\). Since 
magnetic geodesics are uniquely determined by their initial values, 
\(\tilde{\gamma}\) and \(\hat{\gamma}\) coincide on their common domain of 
definition. This yields a continuation of \(\gamma\)
as a magnetic geodesic on \((a,b+\beta)\), which contradicts the 
maximality of \(b\).
\end{proof}

\begin{remark}
We will be looking for closed solutions of \eqref{magnetic-geodesic}, 
which Theorem \ref{thm:1} does not inform us about.
Theorem \ref{thm:1} can be generalized to higher dimensions.
Note that Theorem \ref{thm:1} no longer holds on a surface that is in some 
way not complete, for example a surface with singularities.
\end{remark}

Again, since magnetic geodesics are uniquely determined by their initial values, 
the intermediate value theorem gives us 
the following method for finding closed magnetic geodesics, which 
was employed to produce the numerical examples of closed magnetic 
geodesics found in the figures in this paper: 

\begin{proposition}\label{thm:2}
Let $n=2$.  
Suppose there exists a continuous one-parameter family of solutions 
$\gamma_s$, with $s \in [0,1]$, as in Theorem \ref{thm:1}, and 
suppose there exist $t_1(s)$ and $t_2(s)$ in $\mathbb{R}$ with 
$t_2(s) > t_1(s)$ such that
\begin{enumerate}
\item $t_1(s)$ and $t_2(s)$ depend continuously on $s$, 
\item $\gamma_s(t_1(s)) = \gamma_s(t_2(s))$ for all $ s \in [0,1]$, 
\item $\{ \gamma_0^\prime(t_1(0)) , \gamma_0^\prime(t_2(0))  \}$ spans 
$T_{\gamma_0(t_1(0))}M=T_{\gamma_0(t_2(0))}M$ with one 
orientation, \\
and $\{ \gamma_1^\prime(t_1(1)) , \gamma_1^\prime(t_2(1))  \}$ spans 
$T_{\gamma_1(t_1(1))}M=T_{\gamma_1(t_2(1))}M$ with the 
opposite orientation.  
\end{enumerate}
Then $\gamma : [t_1(s),t_2(s)] \to M$ forms a closed loop for some 
$s \in (0,1)$.  
\end{proposition}

For our numerical studies of \eqref{equation-surface} we need the following

\begin{proposition}\label{prop:2pt4}
Let \(M\subset\mathbb{R}^3\) be a surface. Then equation \eqref{equation-surface} is equivalent to the system
\begin{align}
\label{wayne1} |\gamma'|^2&=c \text{ is constant,}\\
\label{wayne2} \frac{1}{|n|}\langle\gamma'',\gamma'\times n\rangle&=\kappa |\gamma'|^2,
\end{align}
where \(n\) denotes a normal to the surface compatible with $J^{90}$ 
and \(\times\) denotes the cross product in \(\mathbb{R}^3\). 
\end{proposition}
\begin{proof}
The first equation can easily be derived from \eqref{equation-surface} (see 
also Equation \eqref{eqn1pt3}):
\[
\frac{\partial}{\partial t}\frac{1}{2}|\gamma'|^2=\langle\nabla_{\gamma'}\gamma',\gamma'\rangle
=\kappa\langle J^{90}_\gamma(\gamma'),\gamma'\rangle=0.
\]
For the second equation, we consider
\[
\frac{1}{|n|}\langle\gamma'',\gamma'\times n\rangle = 
\langle\gamma'', J^{90}_\gamma(\gamma') 
\rangle = \frac{1}{\kappa} \langle\gamma'', 
\nabla_{\gamma'}\gamma' \rangle = \frac{1}{\kappa} 
|\nabla_{\gamma'}\gamma'|^2 \; . 
\]  Since the magnetic geodesic equation implies 
$|\nabla_{\gamma'}\gamma'|^2=\kappa^2 |\gamma'|^2$, we obtain the second 
equation.  

To establish the equivalence between \eqref{magnetic-geodesic} and the system \eqref{wayne1}, \eqref{wayne2} we note 
that \eqref{wayne1}, \eqref{wayne2} is obtained from \eqref{magnetic-geodesic} by taking the scalar product with
both \(\gamma'\) and \(J_\gamma^{90}(\gamma')\). However, \(\gamma', J_\gamma^{90}(\gamma')\) form a basis of the
tangent space \(T_\gamma M\), yielding the equivalence.
\end{proof}

We now consider a surface $S(u,v)$ parametrized by coordinates $(u,v)$ in a 
subdomain of $\mathbb{R}^2$, and a curve $\gamma(t)=S(u(t),v(t))$ on the surface.  
We can rewrite \eqref{wayne1} and \eqref{wayne2}:  Expanding to obtain 
\begin{align*}
\gamma'&=S_uu'+S_vv',\\
\gamma''&=S_{uu}u'^2+S_{vv}v'^2+S_uu''+S_vv''+2S_{uv}u'v'
\end{align*}
and taking $n=S_u \times S_v$, and using 
\[
\gamma'\times n=\gamma'\times (S_u\times S_v) 
= \langle\gamma',S_v\rangle S_u-\langle\gamma',S_u\rangle S_v \; , 
\]
we can convert equations \eqref{wayne1} and \eqref{wayne2} into 
\begin{align}
\label{volker2}
c=&|S_u|^2u'^2+|S_v|^2v'^2+2\langle S_u,S_v\rangle u'v', \\
\label{volker1}
c |S_u \times S_v| \kappa=&(u''v'-v''u')(|S_v|^2|S_u|^2-|\langle S_u,S_v\rangle|^2)\\
\nonumber &+u'^3(\langle S_u,S_v\rangle\langle S_{uu},S_u\rangle-|S_u|^2\langle S_v,S_{uu}\rangle)
+v'^3(|S_v|^2\langle S_{vv},S_u\rangle -\langle S_v,S_u\rangle\langle S_v,S_{vv}\rangle) \\
\nonumber &+u'^2v'(|S_v|^2\langle S_{uu},S_u\rangle-\langle S_u,S_v\rangle\langle S_{uu},S_v\rangle
+2\langle S_u,S_v\rangle\langle S_{uv},S_u\rangle-2|S_u|^2\langle S_{uv},S_v\rangle) \\
\nonumber &+v'^2u'(\langle S_u,S_v\rangle\langle S_u,S_{vv}\rangle-|S_u|^2\langle S_{vv},S_v\rangle
+2\langle S_u,S_{uv}\rangle |S_v|^2-2\langle S_u,S_v\rangle\langle S_{uv},S_v\rangle\rangle).
\end{align}

However, if the surface is conformally parametrized, that is
\[
\langle S_u,S_v\rangle=0,\qquad |S_u|^2=|S_v|^2=f(u,v) \geq 0 ,
\]
the system \eqref{volker2} and \eqref{volker1} simplifies to 
\begin{align}
\label{eqn:2pt6}
c=&(u'^2+v'^2)f, \\
\label{eqn:2pt7}
c \kappa=&(u''v'-v''u')f-u'^3\langle S_v,S_{uu}\rangle+v'^3\langle S_u,S_{vv}\rangle
-\frac{1}{2}u'^2v'f_u+\frac{1}{2}v'^2u'f_v.
\end{align}

Using the formulations \eqref{volker2}, 
\eqref{volker1}, \eqref{eqn:2pt6} and \eqref{eqn:2pt7}, 
we now use the idea in Proposition \ref{thm:2} to numerically produce 
examples of closed magnetic geodesics.  

\subsection{Example: round sphere}\label{sec:sphere}
\begin{figure}[phbt]\label{fig:sphere}
\begin{center}
\includegraphics[width=0.3\linewidth]{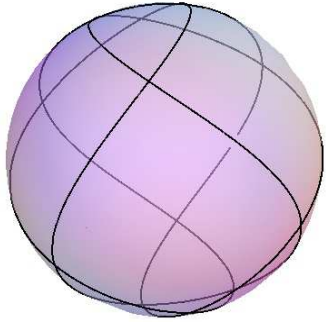}
\caption{A curve with geodesic 
curvature proportional to $\sin u$, on a sphere 
as parametrized in Section \ref{sec:sphere}.}
\end{center}
\end{figure}

Parameterizing the sphere as
\[
S(u,v) = (\cos u \cos v, \cos u \sin v, \sin u) \; , 
\]
the magnetic geodesic system becomes 
\begin{align*}
c=&u'^2+v'^2\cos^2u \; , \\
c \kappa=&(u''v'-v''u')\cos u+v'^3\cos^2u\sin u+2u'^2v'\sin u \; . 
\end{align*}
Note that $\kappa =0$ will give great circles of course, and clearly 
$\kappa$ a nonzero constant will give a circle in the sphere that is not a great circle.  
$\kappa = \sin u$ can give a curve as in Figure \ref{fig:sphere}.  

\subsection{Example: Clifford torus}\label{sec:CliffTorus}
Parameterizing the Clifford torus by
\[
S(u,v) = ((\sqrt{2} + \cos u) \cos v , (\sqrt{2} + \cos u) \sin v, \sin u)
\]
yields the system
\begin{align*}
c=&u'^2+v'^2(\sqrt{2}+\cos u)^2, \\
c \kappa=&(u''v'-v''u')(\sqrt{2}+\cos u) 
+v'\sin u(v'^2(\sqrt{2}+\cos u)^2+2u'^2).
\end{align*}

\begin{figure}[phbt]
\begin{center}
\includegraphics[width=0.4\linewidth]{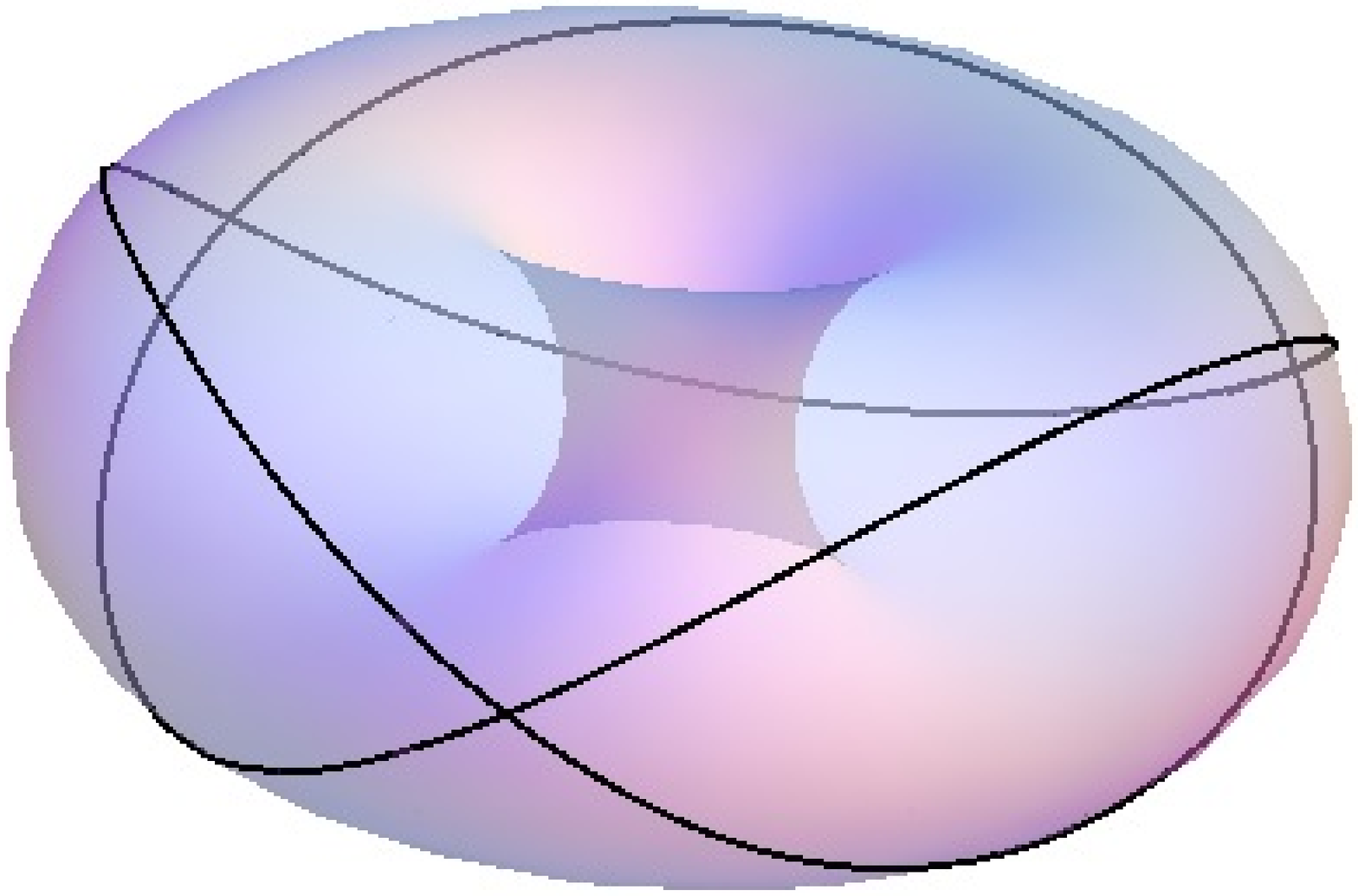}
\includegraphics[width=0.4\linewidth]{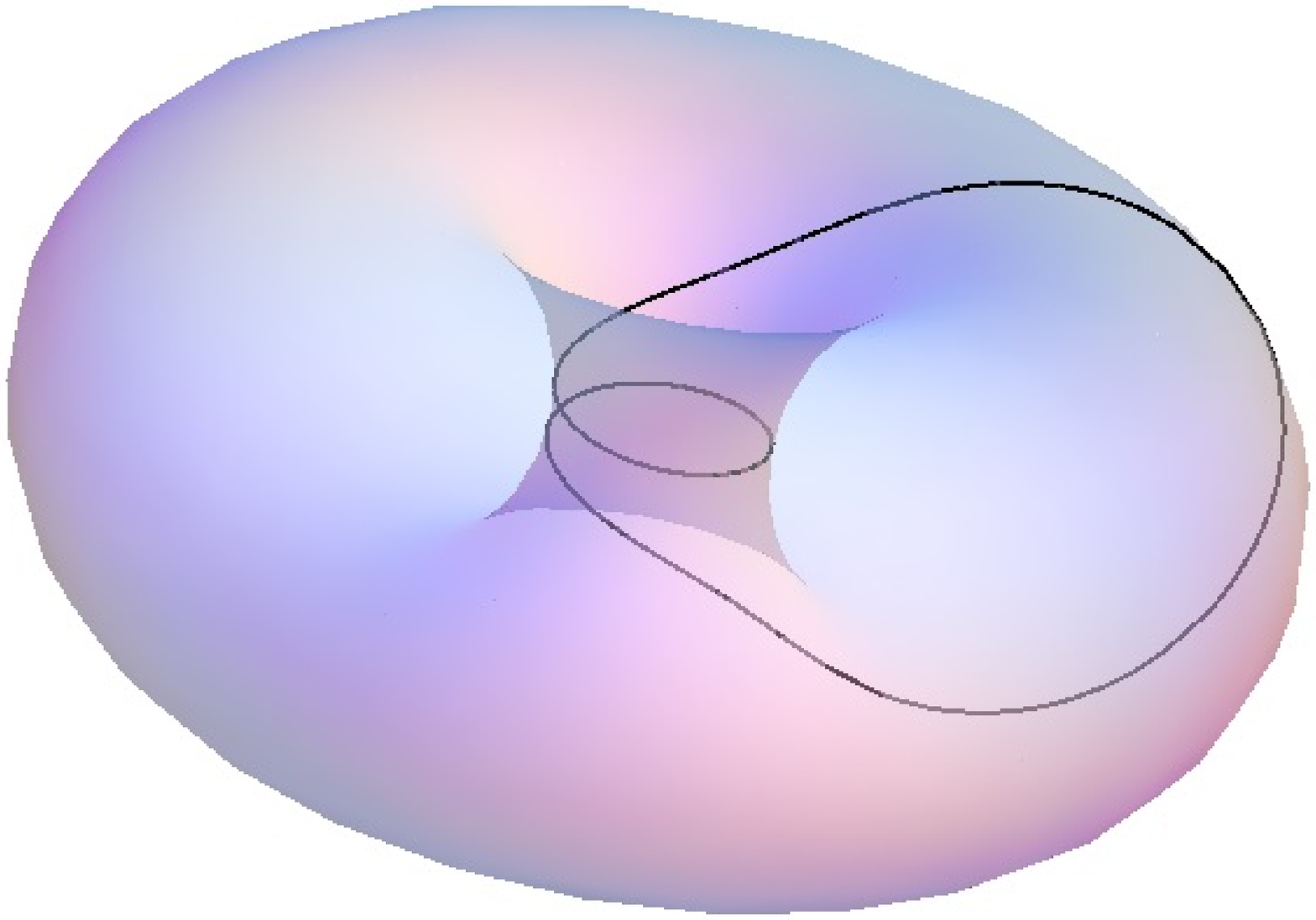}
\caption{Two closed geodesics on the Clifford torus.}
\label{fig:geodonCliff}
\end{center}
\end{figure}

\begin{figure}[phbt]
\includegraphics[width=0.35\linewidth]{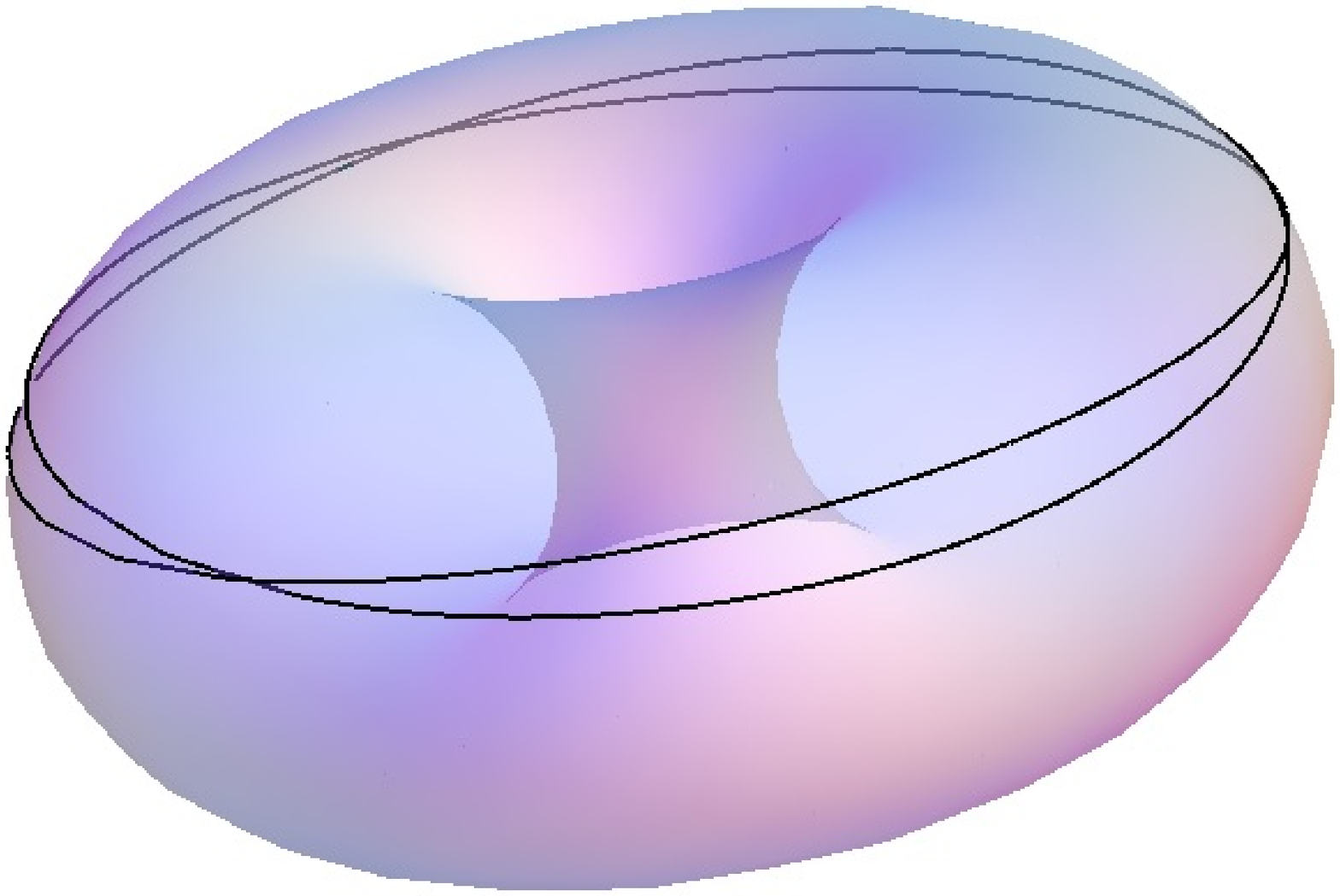}
\includegraphics[width=0.35\linewidth]{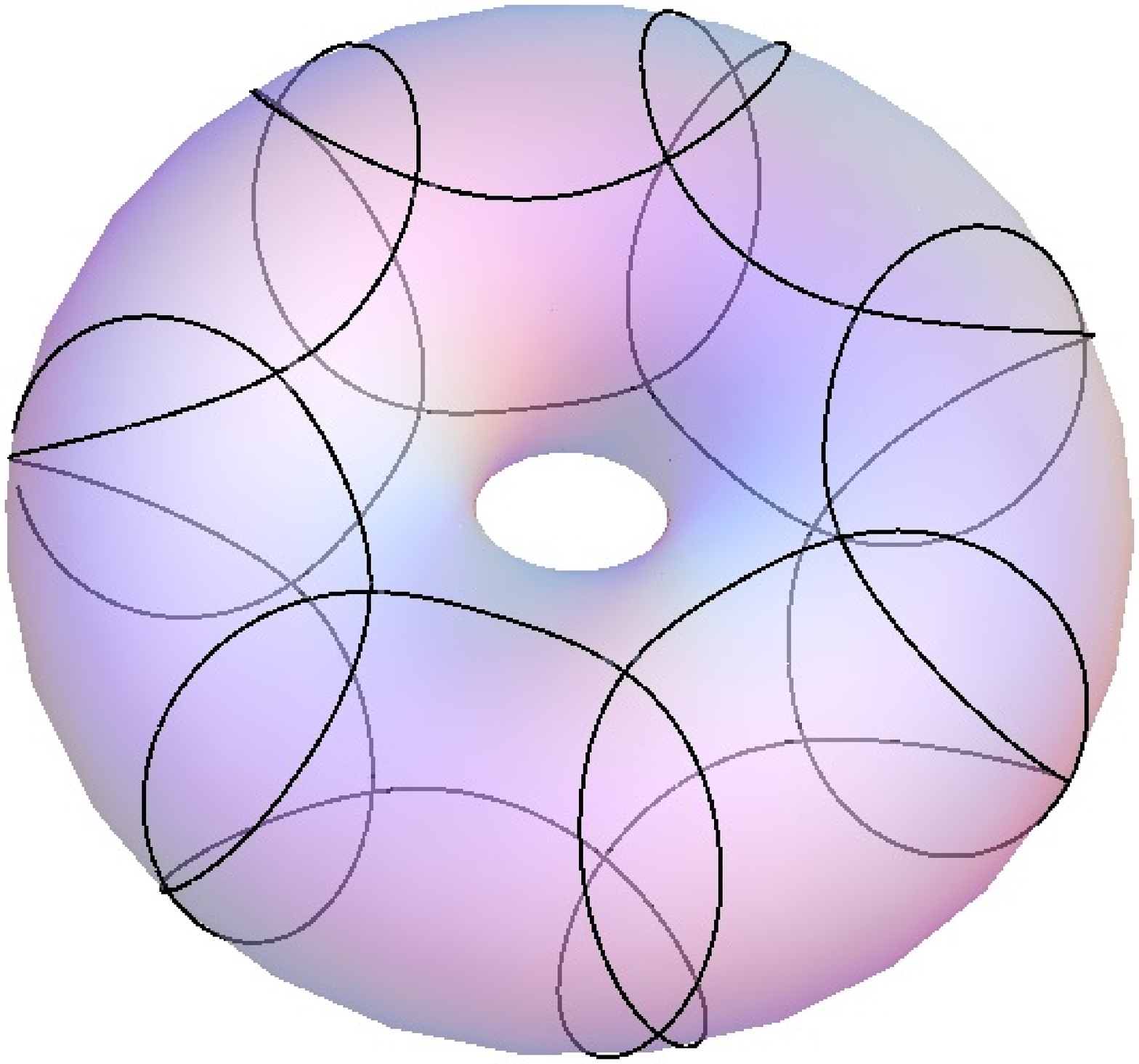}
\includegraphics[width=0.35\linewidth]{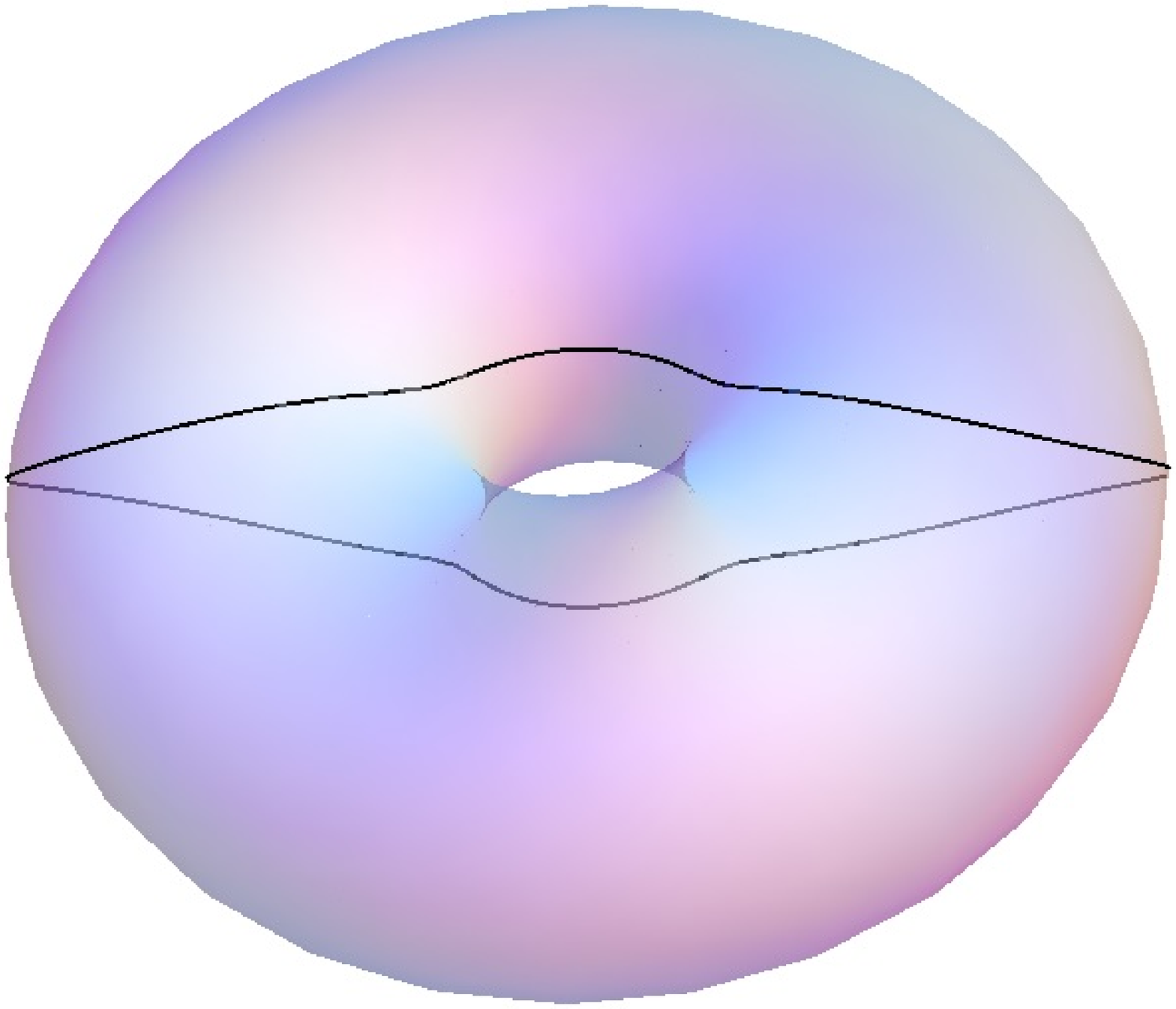}
\caption{
The first picture shows a closed curve with constant non-zero geodesic curvature 
in the Clifford torus, the second picture
a closed curve with geodesic curvature proportional to $\sin u$ in the Clifford torus 
and the third picture another closed curve with geodesic curvature  proportional to 
$\sin u$ in the Clifford torus.}
\label{closedonCliff}
\end{figure}

Two examples of closed geodesics, that is \(\kappa=0\), are given in Figure 
\ref{fig:geodonCliff}.
Other examples of closed magnetic geodesics on the Clifford torus are shown 
in Figure \ref{closedonCliff}.

\subsection{Example: catenoid}\label{sec:cat}
\begin{figure}[phbt]
\begin{center}
\includegraphics[width=0.4\linewidth]{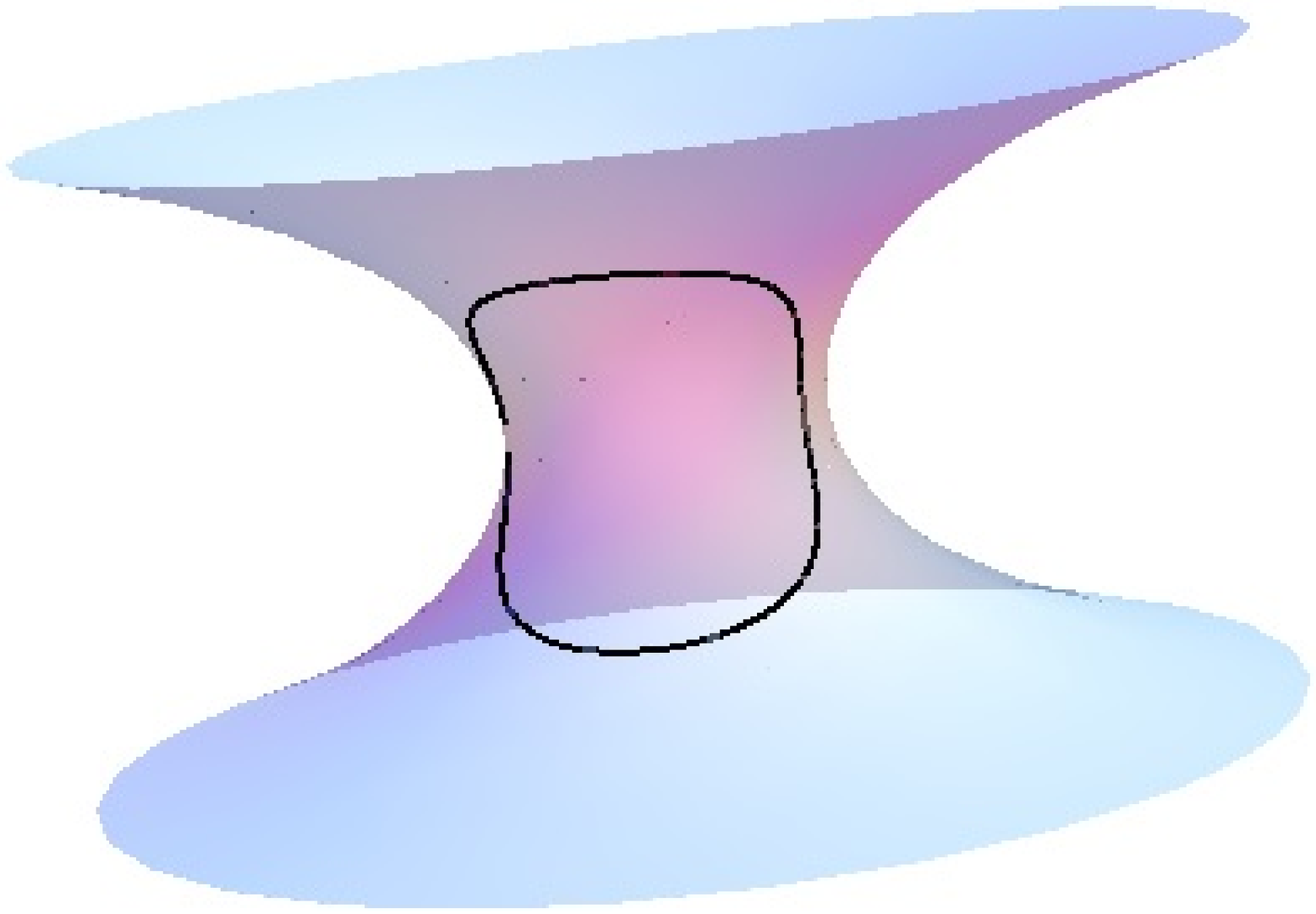}
\includegraphics[width=0.4\linewidth]{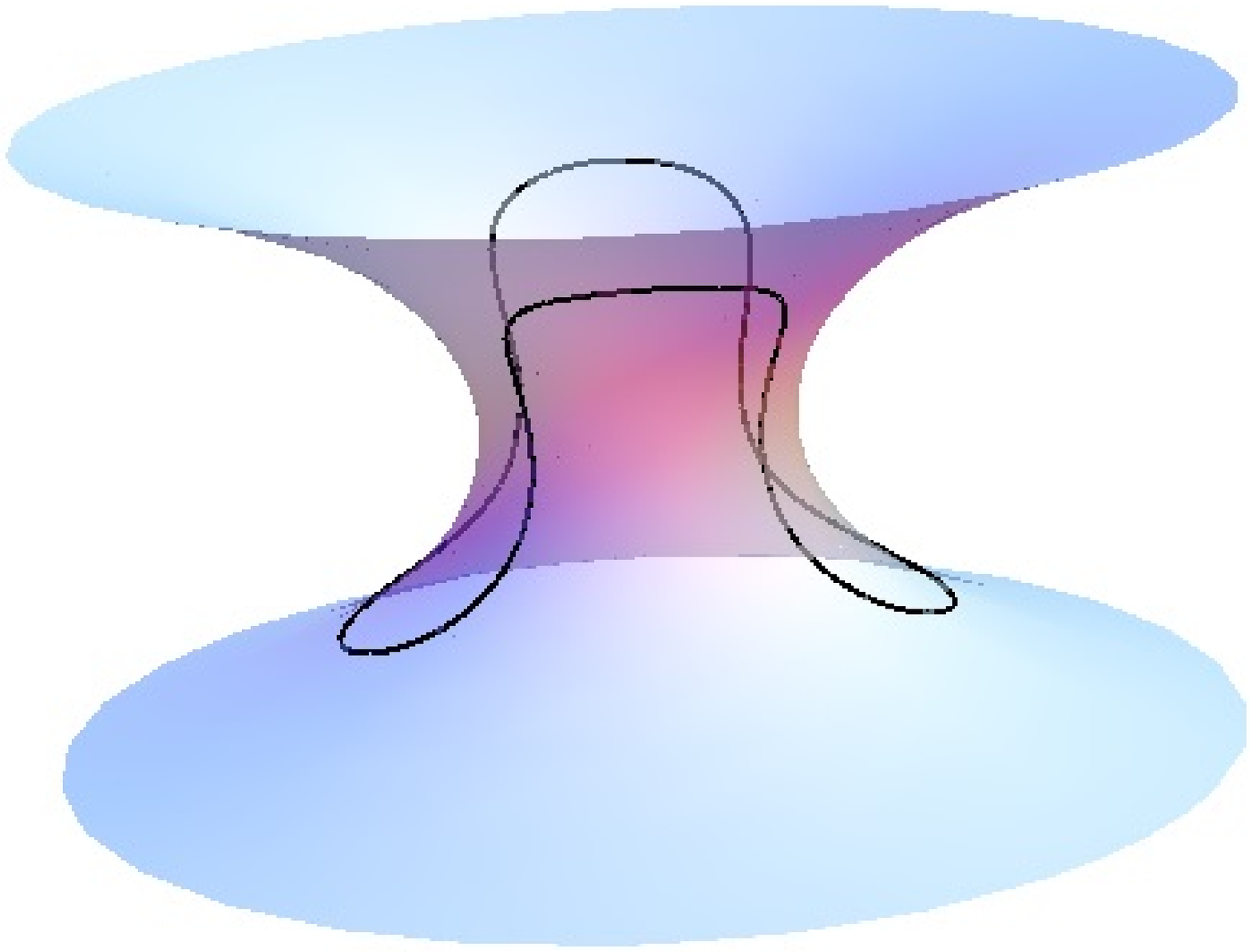}
\caption{The first curve has constant non-zero geodesic curvature on the catenoid,
whereas the second closed curve has geodesic curvature proportional to $\sin u$ on the 
catenoid, as parametrized in Section \ref{sec:cat}.}
\label{fig:cat1}
\end{center}
\end{figure}

Conformally parameterizing the catenoid, with $f=\cosh^2 u$, as 
\[
S(u,v)=(\cosh u \cos v, \cosh u \sin v, u) \; , 
\]
the system becomes 
\begin{align*}
c=&(u'^2+v'^2)\cosh^2u, \\
c \kappa=&(u''v'-v''u')\cosh^2u -
\sinh u\cosh u(v'^3+u'^2v').
\end{align*}
Examples are found in Figure \ref{fig:cat1}.  

\subsection{Example: minimal Enneper surface}
\begin{figure}[phbt]
\begin{center}
\includegraphics[width=0.4\linewidth]{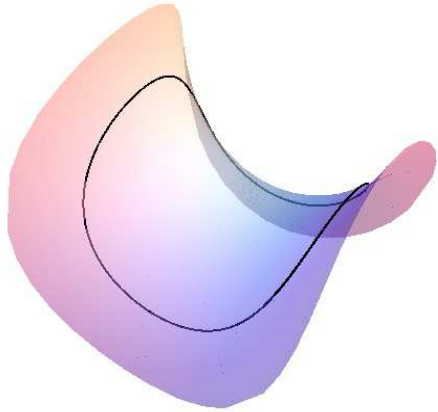}
\caption{A closed curve with constant non-zero geodesic curvature in a 
minimal Enneper surface in $\mathbb{R}^3$.}
\label{fig:EnnExa}
\end{center}
\end{figure}

The Enneper minimal surface in $\mathbb{R}^3$ can be conformally 
parametrized as 
\[
S(u,v)=(u-\frac{1}{3}u^3+u v^2,-v+\frac{1}{3}v^3-v u^2, u^2-v^2) \; , 
\]
with \(f=(1+u^2+v^2)^2\).  This yields the system 
\begin{align*}
c &=(u'^2+v'^2)(1+u^2+v^2)^2, \\
\frac{c \kappa}{1+u^2+v^2}&=(u''v'-v''u')(1+u^2+v^2)+2(u'^2+v'^2)(vu'-uv').
\end{align*}
An example is found in Figure \ref{fig:EnnExa}.  

\subsection{Minkowski \texorpdfstring{$3$}{}-space}
Let $\mathbb{R}^{2,1}$ denote the Minkowski $3$-space $\{  
(x,y,s) \, | \, x,y,s \in \mathbb{R} \}$ with Lorentzian metric of signature $(+,+,-)$.  
Spacelike surfaces with mean curvature identically zero are called \emph{maximal surfaces}, 
and the next example is such a surface.  Our primary result (Theorem \ref{thm:main}) 
is about spacelike surfaces in $\mathbb{R}^{2,1}$, with singularities at 
which the tangent planes become lightlike.  Proposition \ref{prop:2pt4} is true 
for spacelike surfaces in $\mathbb{R}^{2,1}$ as well, once $\mathbb{R}^3$ is 
replaced by $\mathbb{R}^{2,1}$, the cross product for $\mathbb{R}^3 $ is 
replaced by the cross product for $\mathbb{R}^{2,1}$, and the induced connection 
$\nabla$ for surfaces in $\mathbb{R}^3$ is replaced by 
the induced connection $\nabla$ for surfaces in $\mathbb{R}^{2,1}$.  The statement is 
as follows: 

\begin{proposition}\label{prop:2pt4forR21}
Let \(M\subset\mathbb{R}^{2,1}\) be a surface. Then equation \eqref{equation-surface} is equivalent to the system
\begin{align}
|\gamma'|^2&=c \text{ is constant,}\\
\frac{1}{|n|}\langle\gamma'',\gamma'\times n\rangle&=\kappa |\gamma'|^2,
\end{align}
where \(n\) denotes a normal to the surface compatible with $J^{90}$ 
and \(\times\) denotes the cross product in \(\mathbb{R}^{2,1}\). 
\end{proposition}

\subsection{Example: maximal Enneper surface}
In this case we can choose 
\[
S(u,v)=(u+\frac{1}{3}u^3-u v^2, -v-\frac{1}{3}v^3+v u^2, v^2-u^2)
\]
in $\mathbb{R}^{2,1}$.  This parametrization can be obtained from the 
Weierstrass-type representation for maximal surfaces (see, for 
example, \cite{K}), which states that 
\[
S(u,v) = \text{Re}\int^{z=u+iv} (1+g^2,i-ig^2,2g) \eta \; , 
\]
where $g$ is a meromorphic function and $\eta$ is a holomorphic 
$1$-form on a Riemann surface.  This surface is conformally 
parametrized wherever it is nonsingular, and has spacelike 
tangent planes at nonsingular points.  
The singularities occur whenever 
$|g|=1$, and the metric for the surface is 
\[ (1-|g|^2)^2 | \eta |^2 \; . \]
Since, for any magnetic geodesic $\gamma(t)=S(u(t),v(t))$, we have 
\[ c = (u'^2+v'^2) \cdot (1-|g|^2)^2 | \eta |^2 \; , \]
the term $u'^2+v'^2$ 
would have to diverge whenever $\gamma$ approaches 
a singular point.  It follows that magnetic geodesics cannot be extended, 
as solutions of the magnetic geodesic equation, into singular points.  

\begin{figure}[phbt]
\begin{center}
\includegraphics[width=0.35\linewidth]{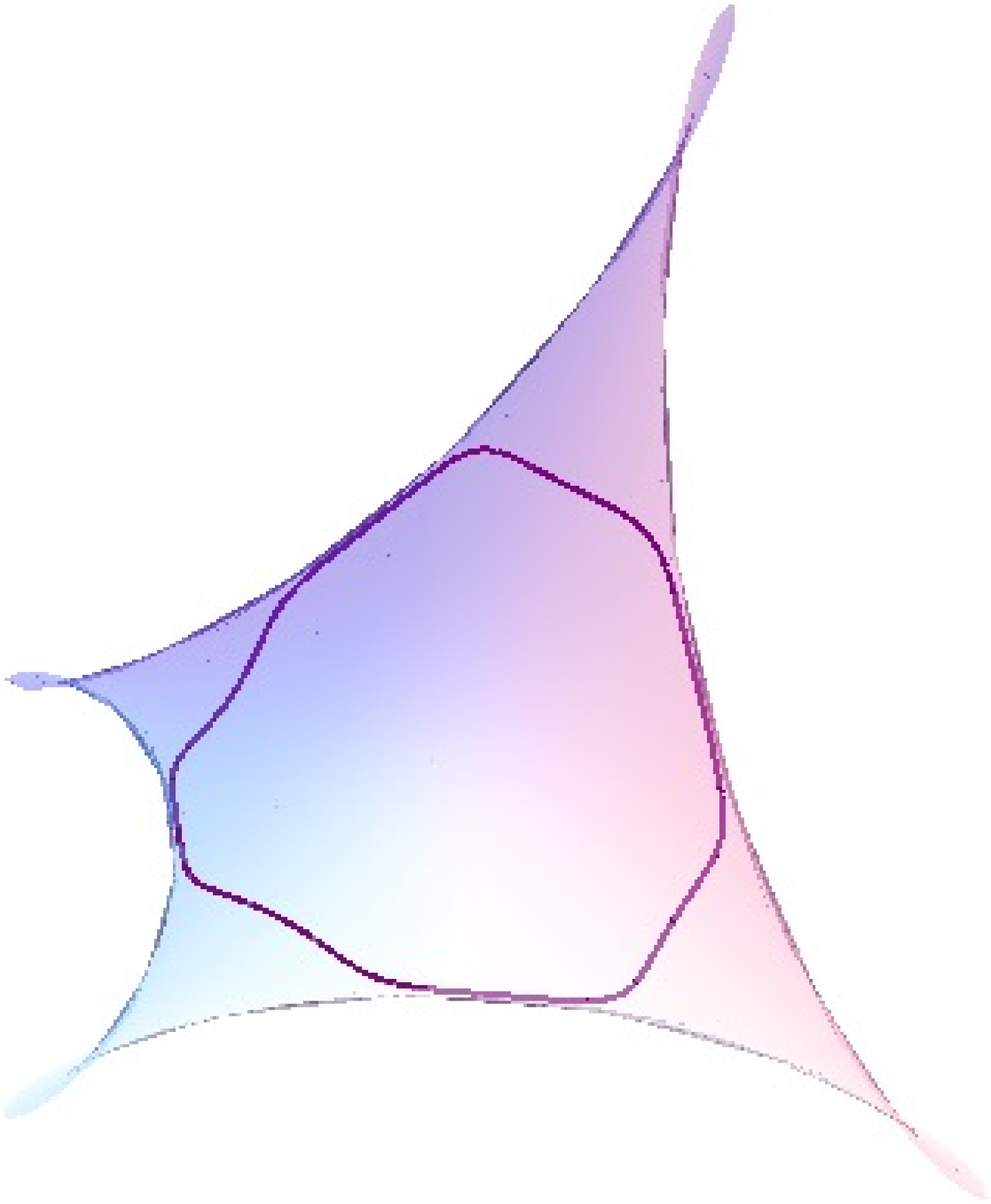}
\includegraphics[width=0.35\linewidth]{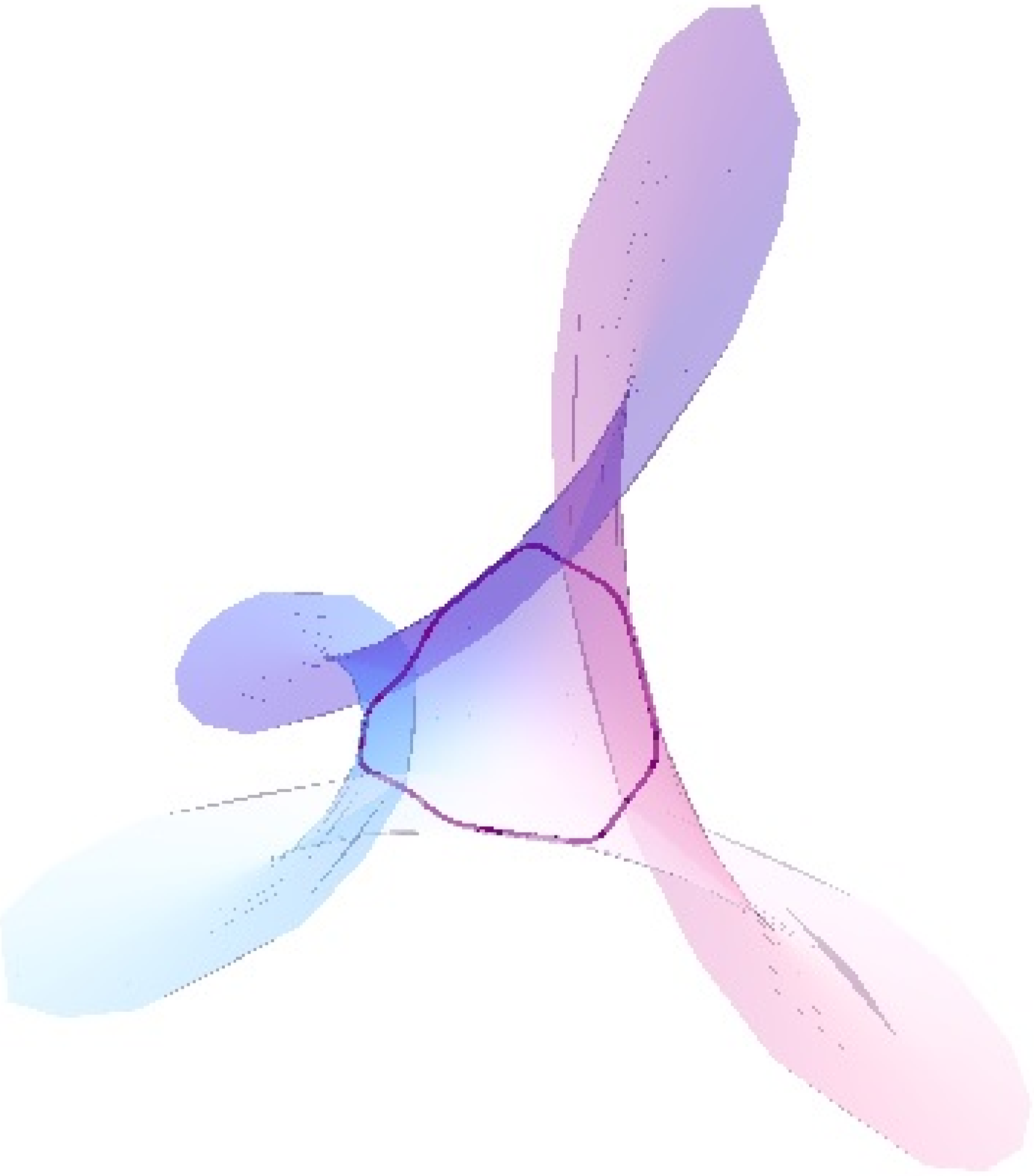}
\caption{A closed curve with constant non-zero geodesic curvature in a 
maximal Enneper surface, shown in both smaller and larger portions of the surface.  
This curve avoids the singular set of the surface.}
\label{fig:maxcat2}
\end{center}
\end{figure}

\begin{figure}[phbt]
\begin{center}
\includegraphics[width=0.35\linewidth]{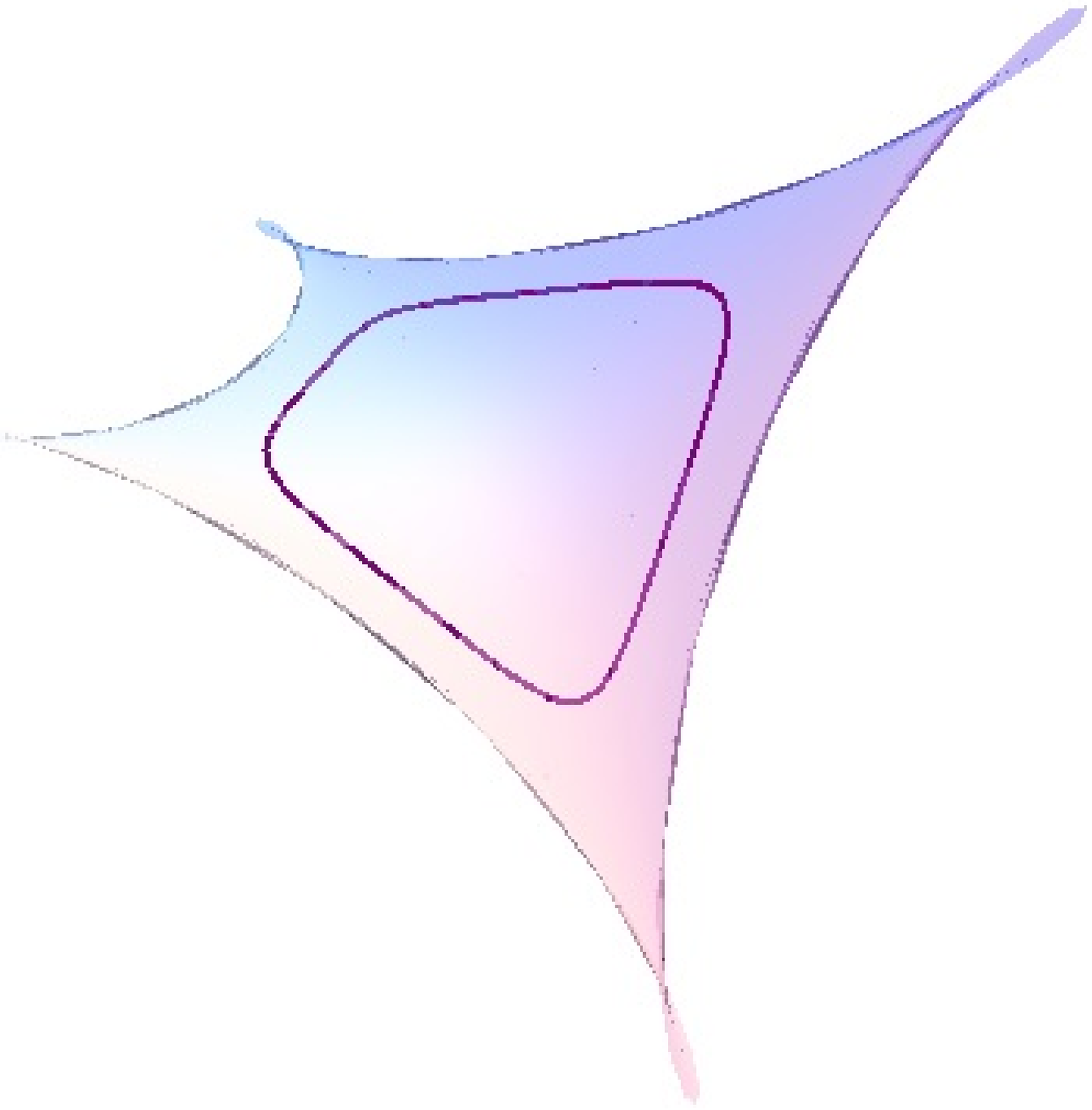}
\includegraphics[width=0.35\linewidth]{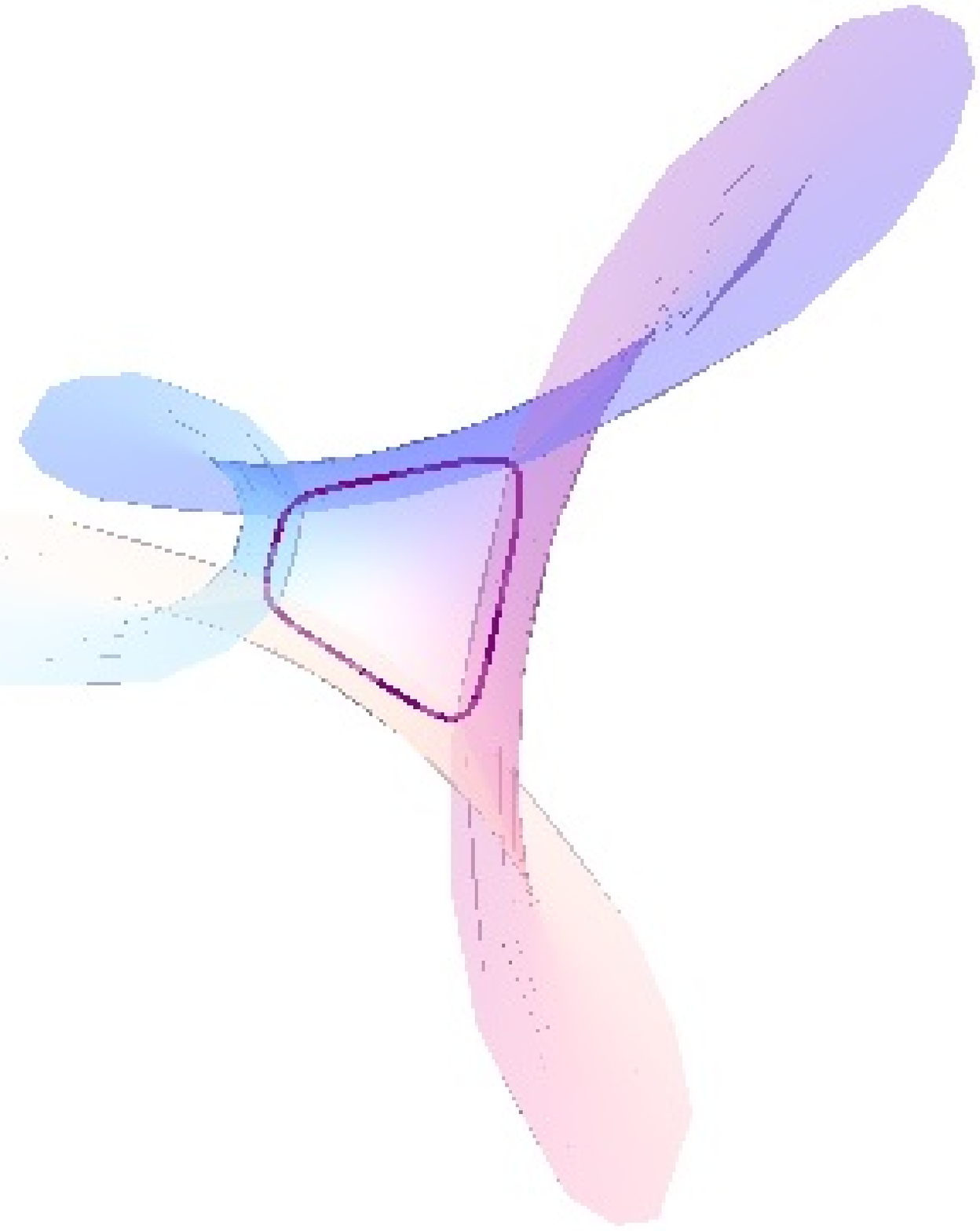}
\caption{A closed geodesic in a maximal Enneper surface in $\mathbb{R}^{2,1}$, 
shown in both smaller and larger portions of that surface.  Note that this geodesic 
also avoids the singular set of the surface.}
\label{fig:maxcat1}
\end{center}
\end{figure}

\begin{figure}[phbt]
\begin{center}
\includegraphics[width=0.35\linewidth]{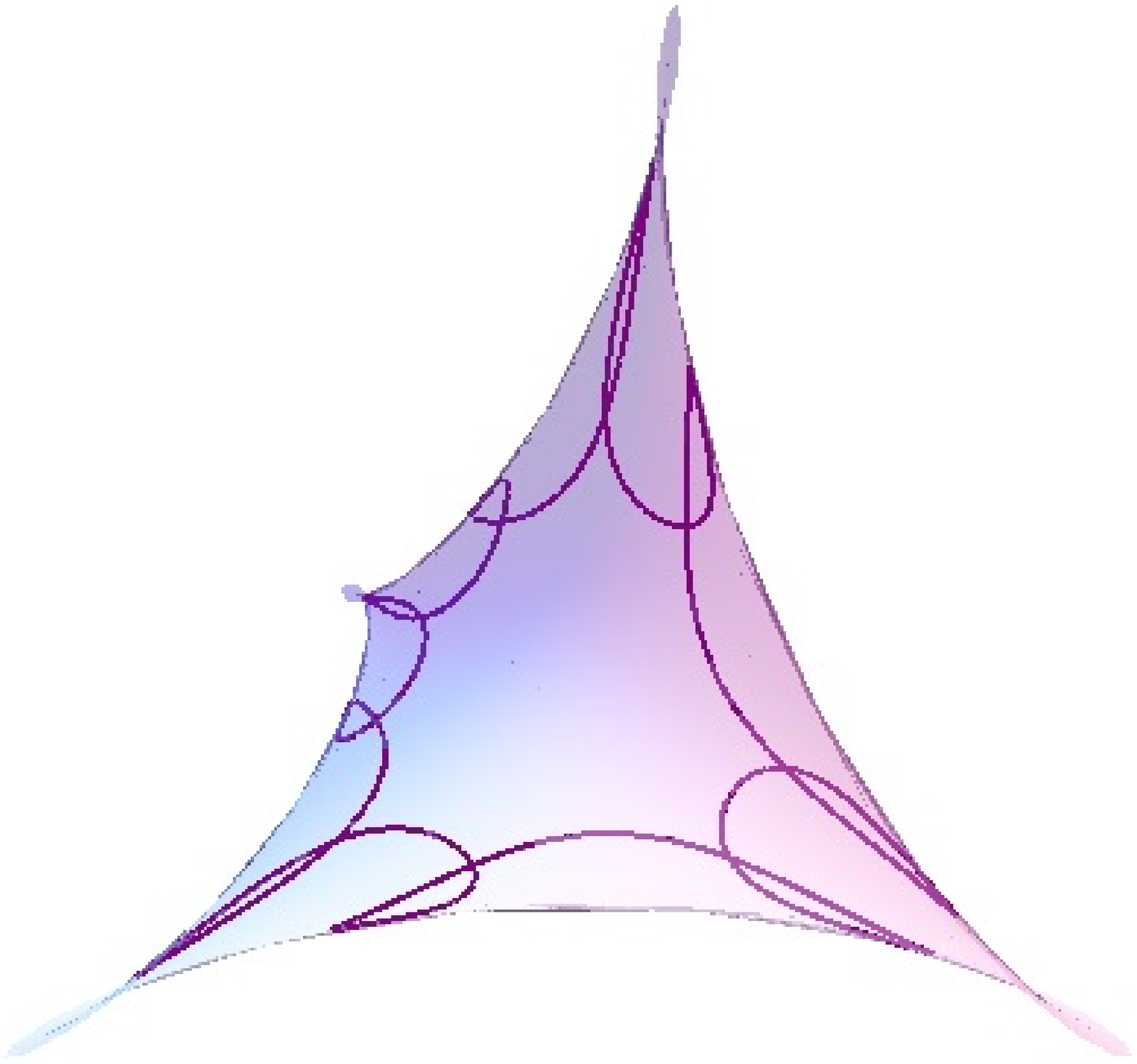}
\includegraphics[width=0.35\linewidth]{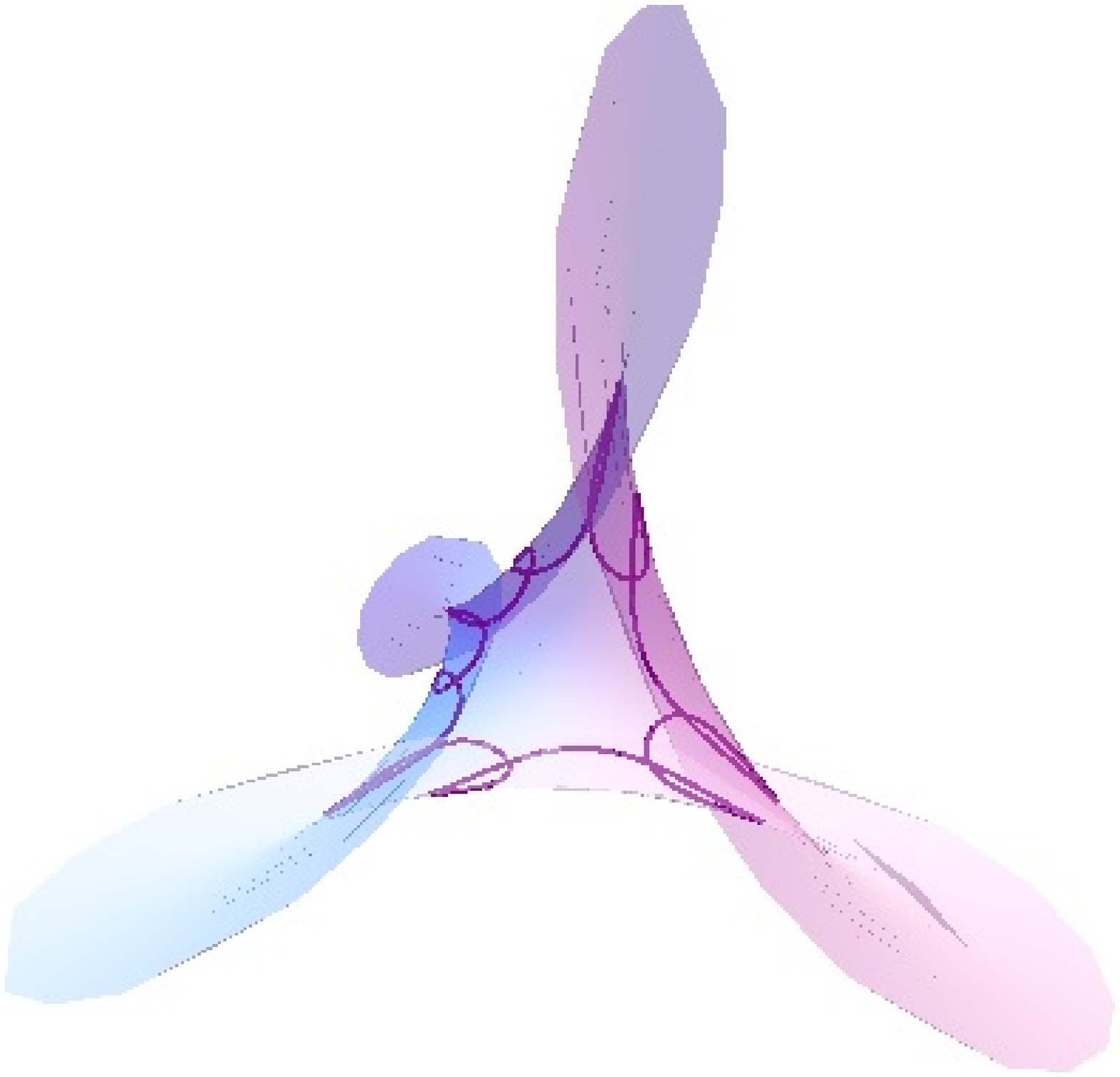}
\caption{Another closed geodesic in a maximal Enneper surface in 
$\mathbb{R}^{2,1}$, again shown in both smaller and larger portions of that 
surface.  Note again that the geodesic avoids the singular set of the surface.}
\label{fig:maxcat3}
\end{center}
\end{figure}

The effect of this fact is 
that magnetic geodesics tend to avoid singular points, as we will see in 
Theorem \ref{thm:main}.  Examples of magnetic geodesics in the maximal Enneper 
surface are shown in Figures \ref{fig:maxcat2}, \ref{fig:maxcat1} and 
\ref{fig:maxcat3}.  

Typically, even at their singularities, maximal surfaces can be described 
as smooth graphs of functions over domains in 
the horizontal spacelike coordinate plane of 
$\mathbb{R}^{2,1}$ (see \cite{FKKRSUYY}, \cite{Gu}, \cite{Kl} for example), 
and thus Theorem \ref{thm:main} will apply to maximal surfaces.  

\subsection{Example: rotated cycloids}
\begin{figure}[phbt]
\begin{center}
\includegraphics[width=0.3\linewidth]{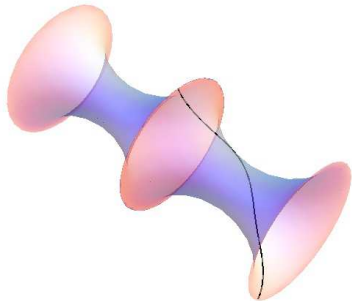}
\caption{A geodesic on a rotated cycloid surface with negative Gaussian curvature.}
\label{fig:cycloid}
\end{center}
\end{figure}

In the case of surfaces in $\mathbb{R}^3$, magnetic geodesics 
will generally not avoid singular sets on those surfaces, and the final example here 
illustrates this.  
We consider rotated cycloids in \(\mathbb{R}^3\), which 
have cuspidal edge singularities.  
We choose the following parametrization
\[
S(u,v)=((2+\cos u) \cos v, (2+\cos u) \sin v, u-\sin u).
\]
The system becomes 
\begin{align*}
c=&2u'^2(1-\cos u)+v'^2(2+\cos u)^2, \\
c \kappa=& (u''v'-v''u')\sqrt{2}\sqrt{1-\cos u}(2+\cos u)+v'^3\frac{\sin u(2+\cos u)^2}{\sqrt{2}\sqrt{1-\cos u}} \\
&+u'^2v'\frac{(6-3\cos u)\sin u}{\sqrt{2}\sqrt{1-\cos u}}.
\end{align*}
An example of a magnetic geodesic that meets the singular set is shown in 
Figure \ref{fig:cycloid}.  

\section{Avoidance of lightlike singularities by magnetic geodesics on surfaces}
In our numerical investigations of magnetic geodesics 
on the maximal Enneper surface we have
seen that magnetic geodesics avoid the singular set of the surface.
In this section we will generalize this conclusion 
not only to arbitrary maximal surfaces, but we 
will also mathematically confirm this behavior on general spacelike surfaces in 
$\mathbb{R}^{2,1}$ at points where the tangent planes degenerate to become lightlike.  
More precisely, we will consider the case 	that the tangent plane $T_pM$ 
becomes lightlike and the surface is a graph of a function over a domain $\mathcal{U}$ 
with immersable boundary $\partial \mathcal{U}$ in the horizontal spacelike 
coordinate plane of $\mathbb{R}^{2,1}$ whose second derivatives are finite and 
not all zero at the projection of $p$ into $\overline{\mathcal{U}}$.  

This is the content of the following theorem: 

\begin{theorem}\label{thm:main}
Suppose that \((M,g)\) is an almost-everywhere-spacelike 
smooth surface in \(\mathbb{R}^{2,1}\) that becomes 
singular at a non-flat point $p \in M$.

Then there are only at most six directions within $T_pM$ to which 
any magnetic geodesic meeting $p$ with $C^1$ regularity and bounded 
geodesic curvature must be tangent.  
Two of these at most six directions are the lightlike directions.  
\end{theorem}

\begin{proof}
We may parametrize the surface as a graph, that is \(S(u,v)=(u,v,f(u,v))\)
for some function \(f(u,v)\), and we can consider a curve $\gamma(t) = S(u(t),v(t))$.  
The surface is spacelike, with the exception of a measure zero set in the surface 
at which the tangent planes are lightlike.  Without loss of generality, we assume 
\begin{enumerate}
\item the tangent plane at \(u=v=0\) is lightlike, 
\item the surface is placed in $\mathbb{R}^{2,1}$ in such a way that 
\[
f(0,0)=0,\qquad f_u(0,0)=1,\qquad f_v(0,0)=0, 
\]
\item the curve $\gamma(t)$ on the surface satisfies 
\[
\gamma(0)=S(0,0),\qquad u'(0)=\cos\theta,\qquad v'(0)=\sin\theta 
\]
for some value of $\theta \in \mathbb{R} \setminus \pi \mathbb{Z}$, 
\item 
the tangent planes to $f$ at the points $\gamma(t)$ for $t > 0$ are spacelike.
\end{enumerate}
We assume that \(\gamma\) is a magnetic geodesic, thus $\langle \gamma' , \gamma' \rangle$ is 
a positive constant for $t>0$. We set
\[ h=(1-f_u^2-f_v^2)^{-1} \; , \;\;\; R=f_{uu}u'^2+2f_{uv}u'v'+f_{vv}v'^2 \; . \]

First, we examine the limiting behavior of $u''(t)$ and $v''(t)$ 
as t approaches \(0\).  
Because $\langle \gamma' , \gamma' \rangle$ is constant 
for $t>0$, by property (3) above we have $\langle \gamma' , \gamma' 
\rangle = \sin^2 \theta$ for all $ t \geq 0$.  
We can assume $|n|=1$ for $t>0$.  We then have 
\begin{equation}\label{doubleprime1}
\langle \gamma'' , \gamma' \rangle = 0 
\end{equation}
and, by Proposition \ref{prop:2pt4forR21}, 
\begin{equation}\label{doubleprime2}
\langle \gamma'' , \gamma' \times n \rangle = \kappa \sin^2 \theta \; . 
\end{equation}
Writing 
\[ \gamma(t) = (u(t),v(t),f(u(t),v(t))) \] 
and using 
\[ \gamma' = (u',v',f_u u'+f_v v') \; , \;\;\; 
\gamma'' = (u'',v'',R+f_u u''+f_v v'') \; , \] we can take 
the limit as $t \to 0$ in Equation \ref{doubleprime1} to obtain the finite limit 
\[
\lim_{t \to 0} (A u'' + B v'')=\cos \theta \cdot R|_{t=0} \; , 
\;\;\;\; A = (1-f_u^2) u' - f_uf_v v' \; , B= (1-f_v^2) v' - f_uf_v u' \; . 
\]
Noting that $A|_{t = 0}=0$ and $B|_{t=0}=\sin \theta \neq 0$, 
we see that only these two cases can occur:
\begin{enumerate}
\item $u''$ is bounded at $t=0$ and $\lim_{t \to 0} v'' = \cot \theta \cdot R|_{t=0}$, or 
\item there exists a sequence $t_j>0$ converging to zero so that 
$|u''(t_j)|$ diverges to infinity and $|v''(t_j)/u''(t_j)|$ converges to zero as $j \to \infty$.   
\end{enumerate}
In the second case, we can obtain the conclusion by examining 
\[ u''(t_j) \left( A|_{t=t_j} + (B|_{t=t_j})\frac{v''(t_j)}{u''(t_j)} \right) \] 
as $j \to \infty$.  

Note that \[ n = \sqrt{h} (f_u,f_v,1) \; , \] so \[ \gamma' \times n = \sqrt{h} 
(f_v (f_u u'+f_v v') -v',u'-f_u (f_u u'+f_v v'),u'f_v-v'f_u) \; . \]  Examining the behavior 
as $t \to 0$ of Equation 
\ref{doubleprime2}, we see that 
\[
{\mathcal T}:= \sqrt{h^{-1}} (u' v''-v' u'') + \sqrt{h} R (v' f_u-u' f_v) \]
is bounded near $t=0$.  

In the first case (1) above with bounded $u''$, ${\mathcal T}$ converges 
asymptotically to $\sqrt{h} R \sin \theta$, and this can be bounded only 
if $R|_{t=0} = 0$.  

In the second case (2) above with unbounded $u''$, we can write ${\mathcal T}$ 
at $t_j$ as 
\[
(u'' \sqrt{h} ( -h^{-1} (v'-u' (v''/u'')) +(R/u'') (v' f_u-u' f_v)))|_{t=t_j} \; . 
\]
Since $u'(t_j)$ and $v'(t_j)$ are bounded, and $v''(t_j)/u''(t_j)$ and $h^{-1}(t_j)$ 
converge to zero, and since 
$(v' f_u-u' f_v)|_{t=t_j}$ converges to $\sin \theta$, as $j \to \infty$, this term is 
asymptotically equal to $(\sqrt{h} R \sin \theta)|_{t=t_j}$ were $R|_{t=0} \neq 0$,  
and again we conclude ${\mathcal T}$ is bounded only if $R|_{t=0} = 0$.  

Thus, in either case, we must have 
\begin{equation}
(f_{uu}\cos^2\theta+2f_{uv}\cos\theta\sin\theta+f_{vv}\sin^2\theta)\big|_{u=v=0} = 0 \; . 
\end{equation}
If \(f_{vv} \neq0\), resp. \(f_{uu} \neq0\), the angle \(\theta\) must satisfy 
\begin{equation}
\tan\theta=\frac{-f_{uv}\pm\sqrt{f_{uv}^2-f_{uu}f_{vv}}}{f_{vv}}\bigg|_{u=v=0} \; , \;\;\; \text{resp.} \;\;\; 
\cot\theta=\frac{-f_{uv}\pm\sqrt{f_{uv}^2-f_{uu}f_{vv}}}{f_{uu}}\bigg|_{u=v=0} \; . 
\end{equation}
If $f_{uu}=f_{vv}=0$, then $\theta = \pi/2+k \pi$ for some integer $k$.  

Thus there are at most four possible values for the angle \(\theta\in[0,2\pi)\)
in addition to \(\theta=0,\pi\) for which the magnetic geodesic can approach the 
singular point $p$.  
\end{proof}

\begin{remark}
Theorem \ref{thm:main} can be 
generalized to almost-everywhere-spacelike submanifolds of 
general dimensional Minkowski spaces, with the corresponding conclusion 
being that generically the possible directions in which a magnetic geodesic can approach a 
point with a lightlike tangent space form a subset in the space of all directions that 
has codimension at least $1$.  
\end{remark}

\end{document}